\documentclass{article}
\usepackage{amsmath,amsthm,amsopn,amstext,amscd,amsfonts,amssymb}
\usepackage{natbib}
\usepackage{verbatim}

\def\tr{\mathop{\rm tr}\nolimits}
\def\R{\mathop{\rm Re}\nolimits}

\def\etr{\mathop{\rm etr}\nolimits}

\renewenvironment{abstract}
                 {\vspace{6pt}
                  \begin{center}
                  \begin{minipage}{5in}
                  \centerline{\textbf{Abstract}}
                  \noindent\ignorespaces
                 }
                 {\end{minipage}\end{center}}
\newtheorem{thm}{\textbf{Theorem}}[section]

\newtheorem{lem}{\textbf{Lemma}}[section]
\theoremstyle{definition}
\newtheorem{defn}{\textbf{Definition}}[section]

\title{\huge \textbf{Bimatrix variate generalised beta distributions}}
\author{
  \textbf{Jos\'e A. D\'{\i}az-Garc\'{\i}a} \thanks{Corresponding author\newline
   {\bf Key words.}  Random matrices, beta distribution, bimatrix variate generalised beta.\newline
    2000 Mathematical Subject Classification. 62E15, 15A52}\\
  Department of Statistics and Computation \\
  25350 Buenavista, Saltillo, Coahuila, Mexico \\
  E-mail: jadiaz@uaaan.mx \\[2ex]
  \textbf{Ram\'on Guti\'errez J\'aimez} \\
  Department of Statistics and O.R. \\
  University of Granada \\
  Granada 18071, Spain \\
  E-mail: rgjaimez@ugr.es\\
}
\date{}
\begin{document}
\maketitle
\begin{abstract}
In this paper, we extend the study of bivariate generalised beta type I and II
  distributions to the matrix variate case.
\end{abstract}

\section{Introduction}

Matrix variate beta type I and II distributions have been studied by different
authors utilising diverse approaches, see \citet{or:64}, \citet{k:70}, \citet{mh:82},
\citet{c:96}, \citet{gn:00}, \citet{dggj:07, dggj:08a, dggj:08b}, among many others.
These distributions play a very important role in various approaches to proving
hypotheses in the context of multivariate analysis, including canonical correlation
analysis, the general linear hypothesis in MANOVA and multiple matrix variate
correlation analysis, see \citet{mh:82} and \citet{sk:79}. All these techniques are
based on the hypothesis that some matrices $\mathbf{A}$ and $\mathbf{B}$ are
independent with Wishart distributions. In the present paper, these results are
generalised, assuming $\mathbf{A}$ and $\mathbf{B}$ to have a matrix variate gamma
distributions.

$\mathbf{A}: m \times m$ is said to have a matrix variate gamma distribution with
parameters $a$ and $m \times m$ positive definite matrix $\mathbf{\Theta}$, this fact
being denoted as $\mathbf{A} \sim \mathcal{G}_{m}(a, \mathbf{\Theta})$, if its
density function is
\begin{equation}\label{gamma}
    \frac{1}{\Gamma_{m}[a]|\mathbf{\Theta}|^{a}}|\mathbf{A}|^{a-(m+1)/2}\etr(-\mathbf{\Theta}^{-1}\mathbf{A})
    (d\mathbf{A}), \quad \mathbf{A} > \mathbf{0},
\end{equation}
where
\begin{equation}\label{measure}
    (d\mathbf{A}) = \bigwedge_{i\leq j}^{m}da_{ij},
\end{equation}
see \citet[pp. 57 and 61]{mh:82} and \citet{gn:00}; where $\Gamma_{m}[a]$ denotes the
multivariate gamma function and is defined as
$$
  \Gamma_{m}[a] = \int_{\mathbf{V}>0} \etr(-\mathbf{V}) |\mathbf{V}|^{a-(m+1)/2} (d\mathbf{V}),
$$
$\R(a) > (m-1)/2$ and $\etr(\cdot) \equiv \exp(\tr(\cdot))$.

As well as the classification of the beta distribution, as beta type I and type II
(see \citet{gn:00} and \citet{sk:79}), two definitions have been proposed for each of
these distributions. Let us focus initially on the beta type I distribution; if
$\mathbf{A}$ and $\mathbf{B}$ have a matrix variate gamma distribution, i.e.
$\mathbf{A} \sim \mathcal{G}_{m}(a,\mathbf{I}_{m})$ and $\mathbf{B} \sim
\mathcal{G}_{m}(b,\mathbf{I}_{m})$ independently, then the beta matrix $\mathbf{U}$
can be defined as
\begin{equation}\label{defbI}
  \mathbf{U} =
  \left\{%
     \begin{array}{ll}
      (\mathbf{A} + \mathbf{B})^{-1/2}\mathbf{A} ((\mathbf{A} + \mathbf{B})^{-1/2})', & \mbox{Definition $1$ or},\\
       \mathbf{A} ^{1/2}(\mathbf{A} + \mathbf{B})^{-1} (\mathbf{A} ^{1/2})', & \mbox{Definition $2$},\\
    \end{array}%
  \right.
\end{equation}
where $\mathbf{C}^{1/2}(\mathbf{C}^{1/2})' = \mathbf{C}$ is a reasonable nonsingular
factorization of $\mathbf{C}$, see \citet{gn:00}, \citet{sk:79} and \citet{mh:82}. It
is readily apparent that under definitions $1$ and $2$ its density function is
denoted as $\mathcal{B}I_{m}(\mathbf{U};a,b)$ and given by
\begin{equation}\label{beta}
     \frac{1}{\beta_{m}[a,b]} |\mathbf{U}|^{(a -(m + 1)/2} |\mathbf{I}_{m} - \mathbf{U}|^{b -(m +1)/2}
     (d\mathbf{U}), \quad \mathbf{0} < \mathbf{U} < \mathbf{I}_{m},
\end{equation}
this being denoted as $\mathbf{U} \sim \mathcal{B}I_{m}(a,b)$, with Re$(a) > (m-1)/2$
and Re$(b)
> (m-1)/2$; where $\beta_{m}[a,b]$ denotes the multivariate beta function defined by
\begin{eqnarray*}
  \beta_{m}[b,a] &=& \int_{\mathbf{0}<\mathbf{S}<\mathbf{I}_{m}} |\mathbf{S}|^{a-(m+1)/2}
  |\mathbf{I}_{m} - \mathbf{S}|^{b-(m+1)/2} (d\mathbf{S}) \\
    &=& \int_{\mathbf{R}>\mathbf{0}} |\mathbf{R}|^{a-(m+1)/2} |\mathbf{I}_{m} + \mathbf{R}|^{-(b+b)} (d\mathbf{R}) \\
    &=& \frac{\Gamma_{m}[a] \Gamma_{m}[b]}{\Gamma_{m}[a+b]}.
\end{eqnarray*}

A similar situation arises with the beta type II distribution, with which we have the
following three definitions:
\begin{equation}\label{defbII}
  \mathbf{F} =
  \left\{%
     \begin{array}{ll}
      \mathbf{B}^{-1/2}\mathbf{A} (\mathbf{B}^{-1/2})', & \mbox{Definition 1},\\
      \mathbf{A}^{1/2}\mathbf{B}^{-1} (\mathbf{A} ^{1/2})', & \mbox{Definition 2},\\
    \end{array}%
  \right.
\end{equation}
with the distribution being denoted as $\mathbf{F} \sim \mathcal{B}II_{m}(a,b)$. In
this case under definition 1 and 2, the density function of $\mathbf{F}$ is denoted
$\mathcal{B}II_{m}(\mathbf{F};a,b)$ by and defined as
\begin{equation}\label{efe}
    \frac{1}{\beta_{m}[a,b]}  |\mathbf{F}|^{a-(m+1)/2}|\mathbf{I}_{m} + \mathbf{F}|^{-(a+b)}
    (d\mathbf{F}), \quad \mathbf{F} > 0.
\end{equation}
$(d\mathbf{F})$ is given in analogous form to (\ref{measure}).

Some of these generalisations from a univariate beta distribution to the matrix
variate case are inappropriate because, in some applications, the researcher is
interested in a vector variate, not in a symmetric matrix, see \citet{ln:82}. In
other words, the researcher is interested in a vector, say, $\mathbf{X} = (x_{1},
\dots , x_{m})'$, such that $x_{i}$ has a marginal beta type I or II distribution for
all $i = 1, \dots,m$. In this respect, \cite{ln:82} and \citet{cn:84} proposed a
multivariate version of the beta type I and II distributions. Let us consider the
following bivariate version, see \citet{ol:03} and \citet{ncg:08}.

Let $X_{0}, X_{1},X_{2}$ be distributed as independent gamma random variates with
parameters $a = a_{0}, a_{1}, a_{2}$, respectively and $\mathbf{\Theta} = 1$ in the
three cases, and define
$$
  U_{1} = \frac{X_{1}}{X_{1}+X_{0}}, \qquad U_{2} = \frac{X_{2}}{X_{2}+X_{0}}.
$$
Clearly, $U_{1}$ and $U_{2}$ each have a beta type I distribution, $U_{1} \sim
\mathcal{B}I_{1}(a_{1},a_{0})$ and $U_{2} \sim \mathcal{B}I_{1}(a_{2},a_{0})$, over
$0 \leq u_{1},u_{2} \leq 1$. However, they are correlated such that $(U_{1},U_{2})'$
has a bivariate generalised beta type I distribution over $0 \leq u_{1},u_{2} \leq
1$. The joint density function of $U_{1}$ and $U_{2}$ is
$$
  \frac{u_{1}^{a_{1}-1} u_{2}^{a_{2}-1} (1 - u_{1})^{a_{2} + a_{0}-1} (1 - u_{2})^{a_{1} + a_{0}-1}}
  {\beta^{*}_{1}[a_{1},a_{2},a_{0}] (1 - u_{1}u_{2})^{a_{1}+ a_{2}+ a_{0}}}, \quad 0 \leq u_{1},u_{2} \leq 1
$$
where
$$
  \beta^{*}_{m}[a,b,c] = \frac{\Gamma_{m}[a] \Gamma_{m}[b]
  \Gamma_{m}[c]}{\Gamma_{m}[a+b+c]}.
$$

A similar result is obtained in the case of beta type II. Here define
$$
  F_{1} = \frac{X_{1}}{X_{0}}, \qquad F_{2} = \frac{X_{2}}{X_{0}}.
$$
Once again it is evident that $F_{1}$ and $F_{2}$ each have a beta type II
distribution, $F_{1} \sim \mathcal{B}II_{1}(a_{1},a_{0})$ and $F_{2} \sim
\mathcal{B}II_{1}(a_{2},a_{0})$, over $f_{1},f_{2} \geq 0$. As in the beta type I
case, they are correlated such that $(F_{1},F_{2})'$ has a bivariate generalised beta
type II distribution over $f_{1},f_{2} \geq 0$. The joint density function of $F_{1}$
and $F_{2}$ is
$$
  \frac{f_{1}^{a_{1}-1} f_{2}^{a_{2}-1}}
  {\beta^{*}_{1}[a_{1},a_{2},a_{0}] (1 + f_{1} + f_{2})^{a_{1}+ a_{2}+ a_{0}}},
  \quad f_{1},f_{2} \geq 0
$$
Some applications to utility modelling  and Bayesian analysis are presented in
\citet{ln:82} and \citet{cn:84}, respectively. Properties such as the moments
$u_{1}^{r}u_{2}^{s}$, conditional distribution, the distributions of the product
$u_{1}u_{2}$, and the $u_{1}/u_{2}$ and $u_{1}/(u_{1} + u_{2})$ quotients are studied
in \citet{ln:82}, \citet{cn:84}, \citet{ol:03} and \citet{ncg:08}.

In the present paper, we extend the bivariate generalised beta type I and II
distributions to the matrix variate case, see Section \ref{sec3} and \ref{sec4}.
These distributions are termed as bimatrix variate generalised beta type I and II
distributions. In Section \ref{sec5}, some properties of these distributions are
studied.

\section{Preliminary results}\label{sec2}

In this section, some results for the hypergeometric function with a matrix argument
are shown.

\begin{defn}\label{def1}
The hypergeometric functions of a matrix argument are given by
\begin{equation}\label{hf}
    {}_{p}F_{q}(a_{1},\dots,a_{p};b_{1},\dots,b_{q};\mathbf{X}) = \sum_{t=0}^{\infty} \sum_{\tau}
    \frac{(a_{1})_{\tau} \cdots (a_{p})_{\tau}}{(b_{1})_{\tau} \cdots (b_{q})_{\tau}}
    \frac{C_{\tau}(\mathbf{X})}{t!}
\end{equation}
where $\sum_{\tau}$ denotes the summation over all the partitions $\tau = (t_{1},
\dots, t_{m})$, $t_{1} \geq \cdots \geq t_{m} \geq 0$, of $t$, $C_{\tau}(\mathbf{X})$
is the zonal polynomial of $\mathbf{X}$ corresponding to $\tau$ and the generalised
hypergeometric coefficient $(a)_{\tau}$ is given by
$$
  (a)_{\tau} = \prod_{i = 1}^{m} (a - (i-1)/2)_{t_{i}},
$$
where $(a)_{t} = a(a+1)(a+2) \cdots (a+t-1)$, $(a)_{0} = 1$. Here $\mathbf{X}$, the
argument of the function, is a complex symmetric $m \times m$ and the parameters
$a_{i}$, $b_{j}$ are arbitrary complex numbers.
\end{defn}
Some other characteristics of the parameters $a_{i}$ and $b_{j}$ and the convergence
of (\ref{hf}) appear in \citet[p. 258]{mh:82}.

A special case of (\ref{hf}) is
\begin{eqnarray*}
  {}_{1}F_{0}(a;\mathbf{X}) &=& \sum_{t=0}^{\infty} \sum_{\tau}
    (a)_{\tau}  \frac{C_{\tau}(\mathbf{X})}{t!}  \qquad (\|\mathbf{X}\| < 1)\\
    &=& |\mathbf{I}_{m} - \mathbf{X}|^{-a}
\end{eqnarray*}
where $\|\mathbf{X}\|$ denotes the maximum of the absolute values of the eigenvalues
of $\mathbf{X}$.

An interesting relation is derived from the following Lemma \ref{lemma}, which gives
an induction method for constructing hypergeometric functions, i.e. integration
involving ${}_{p}F_{q}$ leads to the new hypergeometric function ${}_{p+1}F_{q+1}$. A
motivation for the general recursion comes from the following well-known expressions,
see \citet[Theorem 7.4.2, p. 264]{mh:82}:

\begin{eqnarray*}
    {}_{1}F_{1}(a;c; \mathbf{X})=
    \frac{1}{\beta_{m}[a,c-a]}\hspace{6cm}\\
    \times \int_{0<\mathbf{Y}<\mathbf{I}_{m}} {}_{0}F_{0}(\mathbf{XY})|\mathbf{Y}|^{a-(m+1)/2}
    |\mathbf{I}-\mathbf{Y}|^{c-a-(m+1)/2}(d\mathbf{Y}).
\end{eqnarray*}
and
\begin{eqnarray}
    {}_{2}F_{1}(a,a_{1};c; \mathbf{X})=
    \frac{1}{\beta_{m}[a,c-a]}\hspace{6cm}\nonumber\\ \label{euler}
    \times \int_{0<\mathbf{Y}<\mathbf{I}_{m}} {}_{1}F_{0}(a_{1};\mathbf{XY})|\mathbf{Y}|^{a-(m+1)/2}
    |\mathbf{I}-\mathbf{Y}|^{c-a-(m+1)/2}(d\mathbf{Y}).
\end{eqnarray}
Thus we have
\begin{lem}\label{lemma}
Let $\mathbf{X}<\mathbf{I}$, $\R(a)>(m-1)/2$, $\R(c)>(m-1)/2$ and $\R(c-a)>(m-1)/2$.
Then
\begin{eqnarray*}
    {}_{p+1}F_{q+1}(a,a_{1}, \dots,a_{p};c,b_{1}, \dots, b_{q}; \mathbf{X})=
    \frac{1}{\beta_{m}[a,c-a]}\hspace{4cm}\\
    \times \int_{0<\mathbf{Y}<\mathbf{I}_{m}} {}_{p}F_{q}(a_{1} \cdots a_{p};b_{1} \cdots b_{q};
    \mathbf{XY})|\mathbf{Y}|^{a-(m+1)/2} |\mathbf{I}-\mathbf{Y}|^{c-a-(m+1)/2}(d\mathbf{Y}).
\end{eqnarray*}
\end{lem}
\textit{Proof.} First, we apply an expansion in terms of zonal polynomials
$$
    {}_{p}F_{q}(a_{1}, \dots,a_{p};b_{1}, \dots, b_{q}; \mathbf{XY})=\sum_{t=0}^{\infty}
    \sum_{\tau} \frac{(a_{1})_{\tau} \cdots(a_{p})_{\tau}}{(b_{1})_{\tau} \cdots (b_{q})_{\tau}} \frac{C_{\tau}(\mathbf{XY})}{t!}.
$$
Then, after integrating term by term, see \citet[Theorem 7.2.10, p. 254]{mh:82}, we
have that
\begin{eqnarray*}
    &&\hspace{-1cm}\int_{0<\mathbf{Y}<\mathbf{I}_{m}}{}_{p}F_{q}(a_{1} \cdots a_{p};b_{1} \cdots b_{q};
    \mathbf{XY})|\mathbf{Y}|^{a-(m+1)/2}|\mathbf{I}- \mathbf{Y}|^{c-a-(m+1)/2}(d\mathbf{Y}) \\
    &&\hspace{-0.7cm}=\sum_{t=0}^{\infty}
    \sum_{\tau} \frac{(a_{1})_{\tau} \cdots (a_{p})_{\tau}}{(b_{1})_{\tau} \cdots (b_{q})_{\tau} \ t!}
    \int_{0<\mathbf{Y}<\mathbf{I}_{m}}|\mathbf{Y}|^{a-(m+1)/2}|\mathbf{I}-\mathbf{Y}|^{c-a-(m+1)/2}
    C_{\tau}(\mathbf{XY})(d\mathbf{Y})\\
    &&\hspace{-0.7cm}=  \beta_{m}[a,c-a]\sum_{t=0}^{\infty}\sum_{\tau} \frac{(a)_{\tau}}{(c)_{\tau}}
    \frac{(a_{1})_{\tau} \dots (a_{p})_{\tau}}{(b_{1})_{\tau} \cdots (b_{q})_{\tau} \ t!}
    C_{\tau}(\mathbf{X})\\
    &&\hspace{-0.7cm}= \beta_{m}[a,c-a]\,  {}_{p+1}F_{q+1}(a,a_{1} \cdots a_{p};c,b_{1}
    \cdots b_{q}; \mathbf{X}),
\end{eqnarray*}
and the required result follows. \qed

The use of zonal polynomials and the hypergeometric function with a matrix argument
has only recently been extended; to a large extent this is a derived from the work of
\citet{ke:06}, who in \citet{k:04}, provided a program in MatLab with a very
efficient algorithm for calculating of Jack polynomials (in particular zonal
polynomials) and the hypergeometric function with a matrix argument.

\section{Bimatrix variate generalised beta type I distribution}\label{sec3}

Let $\mathbf{A}$, $\mathbf{B}$ and $\mathbf{C}$ be independent, where $\mathbf{A}
\sim \mathcal{G}_{m}(a, \mathbf{I}_{m})$, $\mathbf{B} \sim \mathcal{G}_{m}(b,
\mathbf{I}_{m})$ and $\mathbf{C} \sim \mathcal{G}_{m}(c, \mathbf{I}_{m})$ with $\R(a)
> (m-1)/2$, $\R(b) > (m-1)/2$ and $\R(c) > (m-1)/2$ and let us define
\begin{equation}\label{bgb1}\hspace{-.7cm}
    \mathbf{U}_{1} = (\mathbf{A}+\mathbf{C})^{-1/2}\mathbf{A}(\mathbf{A} + \mathbf{C})^{-1/2}
  \quad \mbox{and} \quad \mathbf{U}_{2} = (\mathbf{B}+\mathbf{C})^{-1/2}\mathbf{B}
  (\mathbf{B}+\mathbf{C})^{-1/2}
\end{equation}
Of course, $\mathbf{U}_{1} \sim \mathcal{B}I_{m}(a,c)$ and $\mathbf{U}_{2} \sim
\mathcal{B}I_{m}(b,c)$. However, they are correlated such that the distribution of
$\mathbb{U} = (\mathbf{U}_{1}\vdots \mathbf{U}_{2})' \in \Re^{2m \times m}$ can be
termed a bimatrix variate generalised beta type I distribution, denoted as
$\mathbb{U} \sim \mathcal{BGB}I_{2m \times m}(a,b,c)$.

\begin{thm}\label{teob1}
Assume that $\mathbb{U} \sim \mathcal{BGB}I_{2m \times m}(a,b,c)$. Then its density
function is
\begin{equation}\label{density1}\hspace{-0.7cm}
  \frac{|\mathbf{U}_{1}|^{a-(m+1)/2} |\mathbf{U}_{2}|^{b-(m+1)/2} |\mathbf{I}_{m} - \mathbf{U}_{1}|
  ^{b + c -(m+1)/2} |\mathbf{I}_{m} - \mathbf{U}_{2}|^{a + c -(m+1)/2}}{\beta^{*}_{m}[a,b,c]
  |\mathbf{I}_{m} - \mathbf{U}_{1}\mathbf{U}_{2}|^{a + b + c }}(d\mathbb{U})
\end{equation}
$\mathbf{0} < \mathbf{U}_{1}< \mathbf{I}_{m}$, $\mathbf{0} < \mathbf{U}_{2}<
\mathbf{I}_{m}$, where the measure
$$
  (d\mathbb{U}) = (d\mathbf{U}_{1})\wedge (d\mathbf{U}_{2}).
$$
and $\R(a)> (m-1)/2$, $\R(b) > (m-1)/2$ and $\R(c) > (m-1)/2$.
\end{thm}
\textit{Proof.} The joint density of $\mathbf{A}$, $\mathbf{B}$ and $\mathbf{C}$ is
$$
  \frac{|\mathbf{A}|^{a-(m+1)/2} |\mathbf{B}|^{b-(m+1)/2} |\mathbf{C}|^{c-(m+1)/2}}{\Gamma_{m}[a]
  \Gamma_{m}[b] \Gamma_{m}[c]} \etr(-(\mathbf{A} + \mathbf{B} + \mathbf{C})) (d\mathbf{A})
  (d\mathbf{B})(d\mathbf{C})
$$
By effecting the change of variable (\ref{bgb1}), then
$$
  (d\mathbf{A})(d\mathbf{B})(d\mathbf{C}) = |\mathbf{C}|^{m+1} |\mathbf{I}_{m} - \mathbf{U}_{1}|
  ^{-(m+1)} |\mathbf{I}_{m} - \mathbf{U}_{2}|^{-(m+1)}
  (d\mathbf{U}_{1})(d\mathbf{U}_{2})(d\mathbf{C}).
$$
The joint density of $\mathbf{U}_{1}$, $\mathbf{U}_{2}$ and $\mathbf{C}$ is
$$
  \frac{|\mathbf{U}_{1}|^{a-(m+1)/2} |\mathbf{U}_{2}|^{b-(m+1)/2} }{\Gamma_{m}[a]
  \Gamma_{m}[b] \Gamma_{m}[c]|\mathbf{I}_{m} - \mathbf{U}_{1}|
  ^{a + (m+1)/2} |\mathbf{I}_{m} - \mathbf{U}_{2}|^{b +(m+1)/2}}
  |\mathbf{C}|^{a+b+c-(m+1)/2} \hspace{2cm}
$$
$$
  \times \etr\left[-(\mathbf{I}_{m} - \mathbf{U}_{2})^{-1}(\mathbf{I}_{m} - \mathbf{U}_{1}\mathbf{U}_{2})
  (\mathbf{I}_{m} - \mathbf{U}_{1})^{-1}\mathbf{C}\right]
  (d\mathbf{C}) (d\mathbf{U}_{1})(d\mathbf{U}_{2}).
$$
Integrating with respect to $\mathbf{C}$ using
\begin{eqnarray*}
  \int_{\mathbf{C} > \mathbf{0}} |\mathbf{C}|^{a+b+c-(m+1)/2}
  \etr\left[-(\mathbf{I}_{m} - \mathbf{U}_{2})^{-1}(\mathbf{I}_{m} - \mathbf{U}_{1}\mathbf{U}_{2})
  (\mathbf{I}_{m} - \mathbf{U}_{1})^{-1}\mathbf{C}\right] (d\mathbf{C}) \\
    = \Gamma[a+b+c]\frac{|\mathbf{I}_{m} - \mathbf{U}_{1}|
  ^{a + b + c} |\mathbf{I}_{m} - \mathbf{U}_{2}|^{a + b+ c}}{|\mathbf{I}_{m} - \mathbf{U}_{1}\mathbf{U}_{2}|^{a + b + c }}
\end{eqnarray*}
(from (\ref{gamma})) gives the stated marginal density function for $
(\mathbf{U}_{1}\vdots \mathbf{U}_{2})'$. \qed

As in the bivariate case \citep{ol:03}, the joint density (\ref{density1}) can be
represented as a mixture. Let us first note that
\begin{eqnarray*}
     |\mathbf{I}_{m} - \mathbf{U}_{1}\mathbf{U}_{2}|^{-(a + b + c)} &=&  {}_{1}F_{0}(a + b + c
     ;\mathbf{U}_{1}\mathbf{U}_{2})\\
     &=& \sum_{t=0}^{\infty} \sum_{\tau} (a + b + c)_{\tau}
     \frac{C_{\tau}(\mathbf{U}_{1}\mathbf{U}_{2})}{t!}.
\end{eqnarray*}
By substituting in (\ref{density1}) we obtain that the joint density function of
$(\mathbf{U}_{1}\vdots \mathbf{U}_{2})'$ is
$$
   \sum_{t=0}^{\infty} \sum_{\tau} \frac{(a + b + c)_{\tau}}{\beta^{*}_{m}[a,b,c]} |\mathbf{U}_{1}|^{a-(m+1)/2}
   |\mathbf{U}_{2}|^{b-(m+1)/2} |\mathbf{I}_{m} - \mathbf{U}_{1}|
  ^{b + c -(m+1)/2} \phantom{XXXXXXXXXXXXX}
$$
\vspace{-1cm}
$$\hspace{6cm}
  \times |\mathbf{I}_{m} - \mathbf{U}_{2}|^{a + c -(m+1)/2}
     \frac{C_{\tau}(\mathbf{U}_{1}\mathbf{U}_{2})}{t!}.
$$
Moreover
$$
   \sum_{t=0}^{\infty} \sum_{\tau} \frac{(a + b + c)_{\tau}\Gamma_{m}[b+c]\Gamma_{m}[a+c]}{\Gamma_{m}[a+b+c]
   \Gamma_{m}[c]} \mathcal{B}I_{m}(\mathbf{U}_{1};a,b+c) \phantom{XXXXXXXXXXXXXXXXXXXXXX}
$$
\vspace{-1cm}
$$\hspace{8cm}
  \times \  \mathcal{B}I_{m}(\mathbf{U}_{2};b,a+c)\frac{C_{\tau}(\mathbf{U}_{1}\mathbf{U}_{2})}{t!}.
$$

\section{Bimatrix variate generalised beta type II distribution}\label{sec4}

Let $\mathbf{A}$, $\mathbf{B}$ and $\mathbf{C}$ be independent, where $\mathbf{A}
\sim \mathcal{G}_{m}(a, \mathbf{I}_{m})$, $\mathbf{B} \sim \mathcal{G}_{m}(b,
\mathbf{I}_{m})$ and $\mathbf{C} \sim \mathcal{G}_{m}(c, \mathbf{I}_{m})$ with $\R(a)
> (m-1)/2$, $\R(b) > (m-1)/2$ and $\R(c) > (m-1)/2$ and let us define
\begin{equation}\label{bgb2}
    \mathbf{F}_{1} = \mathbf{C}^{-1/2}\mathbf{A}\mathbf{C}^{-1/2}
  \quad \mbox{and} \quad \mathbf{F}_{2} = \mathbf{C}^{-1/2}\mathbf{B}
  \mathbf{C}^{-1/2}
\end{equation}
Clearly, $\mathbf{F}_{1} \sim \mathcal{B}II_{m}(a,c)$ and $\mathbf{F}_{2} \sim
\mathcal{B}II_{m}(b,c)$. But they are correlated and so the distribution of
$\mathbb{F} = (\mathbf{F}_{1}\vdots \mathbf{F}_{2})' \in \Re^{2m \times m}$ can be
termed a bimatrix variate generalised beta type II distribution, which is denoted as
$\mathbb{F} \sim \mathcal{BGB}II_{2m \times m}(a,b,c)$.

\begin{thm}\label{teob2}
Assume that $\mathbb{F} \sim \mathcal{BGB}I_{2m \times m}(a,b,c)$. Then its density
function is
\begin{equation}\label{density2}
    \frac{|\mathbf{F}_{1}|^{a-(m+1)/2} |\mathbf{F}_{2}|^{b-(m+1)/2} }{\beta^{*}_{m}[a,b,c]
  |\mathbf{I}_{m} + \mathbf{F}_{1} + \mathbf{F}_{2}|^{a + b + c }}(d\mathbb{F})
\end{equation}
$\mathbf{F}_{1} > \mathbf{0}$, $\mathbf{F}_{2} > \mathbf{0}$, where the measure
$$
  (d\mathbb{F}) = (d\mathbf{F}_{1})\wedge (d\mathbf{F}_{2}).
$$
and $\R(a)> (m-1)/2$, $\R(b) > (m-1)/2$ and $\R(c) > (m-1)/2$.
\end{thm}
\textit{Proof.} As an alternative to proceeding as in Theorem \ref{teob1}. Let us
recall that if $\mathbf{U} \sim \mathcal{B}I_{m}(a,b)$, then
$(\mathbf{I}_{m}-\mathbf{U})^{-1} -\mathbf{I}_{m} \sim \mathcal{B}II_{m}(a,b)$, see
\citet{sk:79} and \citet{dggj:07}. Then
$$
  \mathbb{F} =
  \left (
    \begin{array}{c}
      \mathbf{F}_{1} \\
      \mathbf{F}_{2}
    \end{array}
  \right ) =
  \left (
    \begin{array}{c}
      (\mathbf{I}_{m}-\mathbf{U}_{1})^{-1}-\mathbf{I}_{m} \\
      (\mathbf{I}_{m}-\mathbf{U}_{2})^{-1}-\mathbf{I}_{m}
    \end{array}
  \right )
$$
with the Jacobian given by
$$
  (d\mathbf{U}_{1})(d\mathbf{U}_{2}) =  |\mathbf{I}_{m} + \mathbf{F}_{1}|^{-(m+1)}
  |\mathbf{I}_{m} + \mathbf{F}_{2}|^{-(m+1)} (d\mathbf{F}_{1})(d\mathbf{F}_{2}).
$$
Also, note that
\begin{eqnarray*}
  \mathbf{I}_{m}- (\mathbf{I}_{m}+\mathbf{F}_{1})^{-1}  &=& (\mathbf{I}_{m}+\mathbf{F}_{1})^{-1}
  ((\mathbf{I}_{m}+\mathbf{F}_{1}) - \mathbf{I}_{m}) = (\mathbf{I}_{m}+\mathbf{F}_{1})^{-1}\mathbf{F}_{1}\\
  \mathbf{I}_{m}- (\mathbf{I}_{m}+\mathbf{F}_{2})^{-1} &=&
  (\mathbf{I}_{m}+\mathbf{F}_{2})^{-1}\mathbf{F}_{2},
\end{eqnarray*}
Then the joint density of $(\mathbf{F}_{1}\vdots \mathbf{F}_{2})'$ is
$$
  \frac{|\mathbf{F}_{1}|^{a-(m+1)/2} |\mathbf{F}_{2}|^{b-(m+1)/2} |\mathbf{I}_{m}+\mathbf{F}_{1}|^{-(a+b+c)}
  |\mathbf{I}_{m}+\mathbf{F}_{2}|^{-(a+b+c)}}{\beta^{*}_{m}[a,b,c] |\mathbf{I}_{m} - (\mathbf{I}_{m} +
  \mathbf{F}_{1})^{-1}\mathbf{F}_{1} \mathbf{F}_{2}(\mathbf{I}_{m}+\mathbf{F}_{2})^{-1}|^{a + b + c
  }}(d\mathbb{F}).
$$
The desired results is follow noting that
$$
  \frac{|\mathbf{I}_{m}+\mathbf{F}_{1}|^{-1} |\mathbf{I}_{m}+\mathbf{F}_{2}|^{-1}}
  {|\mathbf{I}_{m} - (\mathbf{I}_{m} + \mathbf{F}_{1})^{-1}\mathbf{F}_{1} \mathbf{F}_{2}
  (\mathbf{I}_{m}+\mathbf{F}_{2})^{-1}|} = |\mathbf{I}_{m} + \mathbf{F}_{1} +
  \mathbf{F}_{2}|^{-1}.   \qquad\qquad\mbox{\qed}
$$

Other properties of the distribution $\mathcal{BGB}II_{2m \times m}(a,b,c)$ can be
found in a similar way.

\section{Properties}\label{sec5}

In this section we calculate the moments
$E(|\mathbf{U}_{1}|^{r}|\mathbf{U}_{2}|^{s})$ and the distributions of the product
$\mathbf{Z} = \mathbf{U}_{2}^{1/2} \mathbf{U}_{1} \mathbf{U}_{2}^{1/2}$ and the
inverse $(\mathbf{U}_{1}^{-1}\vdots \mathbf{U}_{2}^{-1})$.

\begin{thm}
Assume that $(\mathbf{U}_{1}\vdots\mathbf{U}_{2}) \sim \mathcal{BGB}I_{2m \times
m}(a,b,c)$ then
$$
  E(|\mathbf{U}_{1}|^{r}|\mathbf{U}_{2}|^{s}) = \frac{\beta_{m}[a+r,b+c] \beta_{m}[b+s,a+c]}{\beta_{m}^{*}[a,b,c]}
  \phantom{XXXXXXXXXXXXXXXXXXXX}
$$
\vspace{-.5cm}
$$\phantom{XXXXXXX}
  \times {}_{3}F_{2}(a+r,b+s,a+b+c;a+b+c+r,a+b+c+s;\mathbf{I}_{m}),
$$
with $\R(b+r) > (m-1)/2$, and $\R(a+c) > (m-1)/2$.
\end{thm}
\textit{Proof.}
$$
   E(|\mathbf{U}_{1}|^{r}|\mathbf{U}_{2}|^{s}) = \frac{1}{\beta_{m}^{*}[a,b,c]} \int_{\mathbf{0}< \mathbf{U}_{1} <
   \mathbf{I}_{m}} |\mathbf{U}_{1}|^{a+r-(m+1)/2} |\mathbf{I}_{m} - \mathbf{U}_{1}|^{b+c-(m+1)/2}
$$
$$
   \times
   \int_{\mathbf{0}< \mathbf{U}_{2} <  \mathbf{I}_{m}} |\mathbf{U}_{2}|^{b+s-(m+1)/2} |\mathbf{I}_{m} -
   \mathbf{U}_{2}|^{a+c-(m+1)/2} |\mathbf{I}_{m} -  \mathbf{U}_{1}
   \mathbf{U}_{2}|^{-(a+b+c)} (d\mathbf{U}_{2}) (d\mathbf{U}_{1}).
$$
From (\ref{euler}) we have
$$
   E(|\mathbf{U}_{1}|^{r}|\mathbf{U}_{2}|^{s}) = \frac{\beta_{m}[a+c,,b+s]}{\beta_{m}^{*}[a,b,c]} \int_{\mathbf{0}< \mathbf{U}_{1} <
   \mathbf{I}_{m}} |\mathbf{U}_{1}|^{a+r-(m+1)/2} |\mathbf{I}_{m} - \mathbf{U}_{1}|^{b+c-(m+1)/2}
$$
$$
  \phantom{XXXXXXXXXXXXXXXX} \times \
   {}_{2}F_{1}(b+s,a+b+c;a+b+c+s;\mathbf{U}_{1}) (d\mathbf{U}_{1}).
$$
Then, integrating using Lemma \ref{lemma} the desired result is obtained. \qed

\begin{thm}
Assume that $(\mathbf{U}_{1}\vdots\mathbf{U}_{2}) \sim \mathcal{BGB}I_{2m \times
m}(a,b,c)$. Then  the density function of $\mathbf{Z} = \mathbf{U}_{2}^{1/2}
\mathbf{U}_{1} \mathbf{U}_{2}^{1/2}$ is
$$
  \frac{\beta_{m}[a+c,b+c] |\mathbf{Z}|^{a-(m+1)/2} |\mathbf{I}_{m}-\mathbf{Z}|^{c-(m+1)/2}}
  {\beta_{m}^{*}[a,b,c]} {}_{2}F_{1}(a+c,a+c;a+b+2c;\mathbf{I}_{m}-\mathbf{Z}) (d\mathbf{Z})
$$
and
$$
  E(|\mathbf{Z}|^{r}) = \frac{\beta_{m}[a+c,b+c] \beta_{m}[a+r,c]}
  {\beta_{m}^{*}[a,b,c]} {}_{3}F_{2}(c,a+c,a+c;a+c+r,a+b+2c;\mathbf{I}_{m})
$$
with $\mathbf{0} < \R(\mathbf{Z}) < \mathbf{I}_{m}$, $\R(a+b) > (m-1)/2$ and $\R(b+c)
> (m-1)/2$.
\end{thm}
\textit{Proof.} Consider the transformation $\mathbf{Z} = \mathbf{U}_{2}^{1/2}
\mathbf{U}_{1} \mathbf{U}_{2}^{1/2}$, $\mathbf{U}_{2} = \mathbf{U}_{2}$, with
$$
  (d\mathbf{U}_{1})(d\mathbf{U}_{2}) = |\mathbf{U}_{2}|^{-(m+1)/2}
  (d\mathbf{Z})(d\mathbf{U}_{2}).
$$
Then the joint density of $\mathbf{Z}$ and $\mathbf{U}_{2}$ is
$$
  \frac{|\mathbf{Z}|^{a-(m+1)/2} |\mathbf{U}_{2}|^{-(a+c)} |\mathbf{I}_{m} - \mathbf{U}_{2}|^{a + c -(m+1)/2}
  |\mathbf{U}_{2}-\mathbf{Z}|^{b + c -(m+1)/2}}{\beta^{*}_{m}[a,b,c] \ \
  |\mathbf{I}_{m} - \mathbf{Z}|^{a + b + c }}(d\mathbf{Z})(d\mathbf{U}_{2}),
$$
$\mathbf{0} < \mathbf{U}_{2} < \mathbf{Z} < \mathbf{I}_{m}$.

Now let us make the change of variable $\mathbf{W} =
(\mathbf{I}_{m}-\mathbf{Z})^{-1/2} (\mathbf{I}_{m}-\mathbf{U}_{2})
(\mathbf{I}_{m}-\mathbf{Z})^{-1/2}$ and $\mathbf{Z} = \mathbf{Z}$. Noting that
$$
  (d\mathbf{Z})(d\mathbf{U}_{2}) = |\mathbf{I}_{m} - \mathbf{Z}|^{(m+1)/2}
  (d\mathbf{Z})(d\mathbf{W}).
$$
Then the joint density of $\mathbf{Z}$ and $\mathbf{W}$ is
$$
  \frac{|\mathbf{Z}|^{a-(m+1)/2} |\mathbf{I}_{m} - \mathbf{Z}|^{c-(m+1)/2}|\mathbf{W}|^{a + c -(m+1)/2}
  |\mathbf{I}_{m}-\mathbf{W}|^{b + c -(m+1)/2}}{\beta^{*}_{m}[a,b,c] \ \ |\mathbf{I}_{m} - (\mathbf{I}_{m} -
  \mathbf{Z})\mathbf{W}|^{a+c} }(d\mathbf{W})(d\mathbf{Z}),
$$
$\mathbf{0} < \mathbf{Z} < \mathbf{I}_{m}$, $\mathbf{0} < \mathbf{W} <
\mathbf{I}_{m}$.

Integrating with respect to $\mathbf{W}$ using (\ref{euler}) with $\R(\mathbf{Z}) <
\mathbf{I}_{m}$, $\R(a+c) > (m-1)/2 $, and $\R(b+c) > (m-1)/2$, we then obtain the
stated marginal density function for $\mathbf{Z}$.

Now in order to find the expression for $E(|\mathbf{Z}|^{r})$, let us first observe
that by $p= 2$ and $q =1$ in Lemma \ref{lemma} we obtain
\begin{eqnarray*}
  {}_{3}F_{2}(a,a_{1},a_{2};c,b_{1}; \mathbf{X})=
  \frac{1}{\beta_{m}[a,c-a]}\hspace{6.5cm}\\
  \times \int_{0<\mathbf{Y}<\mathbf{I}_{m}} {}_{2}F_{1}(a_{1}a_{2};b_{1};
  \mathbf{XY})|\mathbf{Y}|^{a-(m+1)/2} |\mathbf{I}_{m}-\mathbf{Y}|^{c-a-(m+1)/2}(d\mathbf{Y}).
\end{eqnarray*}
Now, by making the transformation $\mathbf{W} = \mathbf{I}_{m} - \mathbf{Y}$, with
$(d\mathbf{Y}) = (d\mathbf{W})$. Then
$$
  {}_{3}F_{2}(a,a_{1},a_{2};c,b_{1}; \mathbf{X})= \frac{1}{\beta_{m}[a,c-a]}
  \hspace{7cm}
$$
\begin{equation}\label{ecz}\hspace{-1cm}
    \times \int_{0<\mathbf{W}<\mathbf{I}_{m}} {}_{2}F_{1}(a_{1}a_{2};b_{1};
  \mathbf{X}(\mathbf{I}_{m} - \mathbf{W}))|\mathbf{I}_{m} - \mathbf{W}|^{a-(m+1)/2}
  |\mathbf{W}|^{c-a-(m+1)/2}(d\mathbf{W}).
\end{equation}
Then the expression for $E(|\mathbf{Z}|^{r})$ follows immediately from the density
function of $\mathbf{Z}$ and (\ref{ecz}). \qed

\begin{thm}
Assume that $(\mathbf{U}_{1}\vdots\mathbf{U}_{2})' \sim \mathcal{BGB}I_{2m \times
m}(a,b,c)$. Then  the density function of $\mathbb{V} = (\mathbf{V}_{1}\vdots
\mathbf{V}_{2})'= ( \mathbf{U}_{1}^{-1}\vdots \mathbf{U}_{2}^{-1})'$ is
$$
   \frac{|\mathbf{V}_{1}|^{-a-(m+1)/2} |\mathbf{V}_{2}|^{-b-(m+1)/2} |\mathbf{I}_{m} - \mathbf{V}_{1}^{-1}|
  ^{b + c -(m+1)/2} |\mathbf{I}_{m} - \mathbf{V}_{2}^{-1}|^{a + c -(m+1)/2}}{\beta^{*}_{m}[a,b,c]
  |\mathbf{I}_{m} - (\mathbf{V}_{1}\mathbf{V}_{2})^{-1}|^{a + b + c }}(d\mathbb{V})
$$
$\mathbf{0} < \mathbf{V}_{1}< \mathbf{I}_{m}$, $\mathbf{0} < \mathbf{V}_{2}<
\mathbf{I}_{m}$, where the measure
$$
  (d\mathbb{V}) = (d\mathbf{V}_{1})\wedge (d\mathbf{V}_{2}).
$$
and $\R(a)> (m-1)/2$, $\R(b) > (m-1)/2$ and $\R(c) > (m-1)/2$.
\end{thm}
\textit{Proof.} The desired result follows immediately by making the change of
variable $(\mathbf{V}_{1}\vdots \mathbf{V}_{2})'= ( \mathbf{U}_{1}^{-1}\vdots
\mathbf{U}_{2}^{-1})'$ in the joint density of $(\mathbf{U}_{1}\vdots
\mathbf{U}_{2})'$, taking into account that $(d\mathbf{U}_{1}) (d\mathbf{U}_{2}) =
|\mathbf{V}_{1}|^{-(m-1)} |\mathbf{V}_{2}|^{-(m-1)}(d\mathbf{V}_{1})
(d\mathbf{V}_{2})$. \qed

\section{Conclusions}

With the algorithms proposed by \cite{ke:06} for the calculation of Jack polynomials
and hypergeometric functions with matrix arguments, together with the MatLab
implementation by \citet{k:04}, it is now is feasible to evaluate in a very efficient
way expressions such as density functions and moments, as shown in the proceding
sections, as well as highly complex expressions.  The theory developed in this paper
has not yet been applied, with the exception of the bivariate case. However, its
potential role is apparent, for example, in multidimensional scaling (MDS) in the
following context. Bimatrix variate generalised beta distributions can be used as
distributions of the matrices of similarities (or dissimilarities) for an individual
when these matrices of similarities have been obtained at two different times. Given
the genesis of bimatrix variate generalised beta distributions, such distributions
may allow us, in some sense, to model the learning problem. Statistical approaches to
MDS have been studied assuming independence between times (without learning) in the
univariate case by \citet{r:82} and by \citet{vma:08}.

Another potential use appears in the context of shape theory, specifically in the
approach known as affine shape or configuration densities. This approach is currently
being studied, see \citet{cdggf:08}. It was first proposed by \citet{gm:93}, and at
present only the Cauchy configuration densities, have been explored. The use of beta
configuration densities,  only was proposed by \citet{gm:93},  but no additional
information was provided. Perhaps bimatrix variate generalised beta distributions can
be obtained as configuration densities if we assume that two figures (images) of a
single individual or a single object, obtained at two different times, are not
independent. Currently, both applications are under study.

\section*{Acknowledgments}

This research work was partially supported  by CONACYT-M\'exico, Research Grant No. \
81512 and IDI-Spain, Grants No. FQM2006-2271 and MTM2008-05785. This paper was
written during J. A. D\'{\i}az- Garc\'{\i}a's stay as a visiting professor at the
Department of Statistics and O. R. of the University of Granada, Spain.

\end{document}